\documentclass[12pt]{amsart}

\usepackage{amssymb}
\usepackage{graphicx}

 \usepackage{url}

\usepackage{graphicx}

\theoremstyle{remark}

\theoremstyle{definition}

\numberwithin{equation}{section}
\numberwithin{theorem}{section}

\begin{document}




\date{Preprint May 31, 2012. To appear in {\it Mathematical Gazette}}


\title{Trigonometry of {\it {\bf The Gold-Bug}}}


%
\author{Erik Talvila}
\address{Department of Mathematics \& Statistics\\
University of the Fraser Valley\\
Abbotsford, BC Canada V2S 7M8}
\email{Erik.Talvila@ufv.ca}



%

%

\begin{abstract}
The classic Edgar Allan Poe story {\it The Gold-Bug} involves digging for pirate
treasure.  Locating the digging sites requires some simple trigonometry.
\end{abstract}

\maketitle

\let\oldsqrt\sqrt
\def\sqrt{\mathpalette\DHLhksqrt}
\def\DHLhksqrt#1#2{%
\setbox0=\hbox{$#1\oldsqrt{#2\,}$}\dimen0=\ht0
\advance\dimen0-0.2\ht0
\setbox2=\hbox{\vrule height\ht0 depth -\dimen0}%
{\box0\lower0.4pt\box2}}

\providecommand{\abs}[1]{\lvert#1\rvert}
\providecommand{\norm}[1]{\lVert#1\rVert}
\section{Introduction}
\begin{quotation}
{\em
The lanterns having been lit, we all fell to work with a zeal
worthy a more rational cause; and, as the glare fell upon our persons and 
implements, I could not help thinking how picturesque a group we composed, 
and how strange and suspicious our labors must have appeared to any 
interloper who, by chance, might have stumbled upon our whereabouts.

We dug very steadily for two hours.}
\end{quotation}

This is an excerpt from the Edgar Allan Poe story {\it The Gold-Bug}
\cite{poe}.
Published in $1843$, it is a gripping tale about Captain Kidd's 
lost pirate 
treasure map and digging for booty on the Carolina shore.  
Locating the digging spots requires 
a little trigonometry.  Poe has our heroes dig two holes but we will
show that they must really overlap.

The main character of the story is William Legrand, who is a nobleman in
impoverished circumstances.  He lives on tiny, mostly unpopulated, 
Charleston Island off the coast of
South Carolina.  With him is the old, freed, slave Jupiter.  It is widely
believed that Legrand is slightly unbalanced and Jupiter takes it upon himself 
to look after him.  The third character is the narrator, an old friend, who
accompanies them on the 
treasure hunt.

Legrand has found a parchment half buried in the sand beside an old ship's boat.  It is written with peculiar symbols and
turns out to be a cryptogram.  He translates it and finds it gives veiled
directions to a secluded spot in the woods on the mainland.  (Poe gives a
lovely explanation of how to crack a cryptogram.)
Here the three
companions dig by lantern light for what they hope is treasure.
 
The gold-bug itself is a living beetle that Legrand has found.  It seems to be
made of real gold.  This gives the story a chilling air of mysticism.

\section{The digging directions}
\begin{quotation}
{\em\ldots main branch seventh limb east side shoot from the left eye of the 
death's-head a bee line from the tree through the shot fifty feet out\ldots}
\end{quotation}

At the digging site in the woods there is an enormous old tulip tree.  
Legrand has Jupiter
climb it.  Going up $70$ feet, he climbs onto the seventh limb up the east 
side of the tree.  Out on the limb a human skull has been fastened by the
pirates.  He drops
the gold-bug through the right eye socket and it falls to the ground.  They then
draw a line a further $50$ feet out from the tree in this direction and dig a 
hole $4$ feet in
diameter and $5$ feet deep.  Alas, no treasure.  They then realise Jupiter has
dropped the bug through the wrong eye socket.  Legrand moves the mark under
the skull $2\;1/2$ inches, as if it had fallen through the other eye socket.
A line is then extended 
$50$ feet through this new point away from the trunk.  Here they
dig a second time.

\section{Trigonometry}
\begin{quotation}
{\em
I dug eagerly, and now and then caught myself actually looking, with 
something that very much resembled expectation, for the fancied treasure\ldots
[W]e were again interrupted by the violent 
howlings of the dog. 
\ldots [H]e made furious resistance, and, leaping into the hole, tore up the 
mould frantically with his claws. In a few seconds he had uncovered a 
mass of human bones\ldots
}
\end{quotation}

Some simple trigonometry shows that, under reasonable assumptions for the
radius of the tree and the distance the skull is from the trunk, the two
holes must have overlapped, even though the narrator says they were 
several yards apart.  

In Figure~1 we have the set up.  Points labeled in the figure are
$O$ (centre of tree), $a$ (drop point through left eye socket),
$b$ (drop point through right eye socket), $A$ (first digging site),
$B$ (second digging site).  Distances are $r$ (radius of tree trunk),
$L$ (distance of each eye socket from tree trunk),
$r_A$ (radius of first hole), 
$r_B$ (radius of second hole).  All distances are in feet, as they are
in the story.  The angle between $OA$ and $OB$ is $2\theta$.
\begin{figure}[p]
{\includegraphics
[width=6in]
{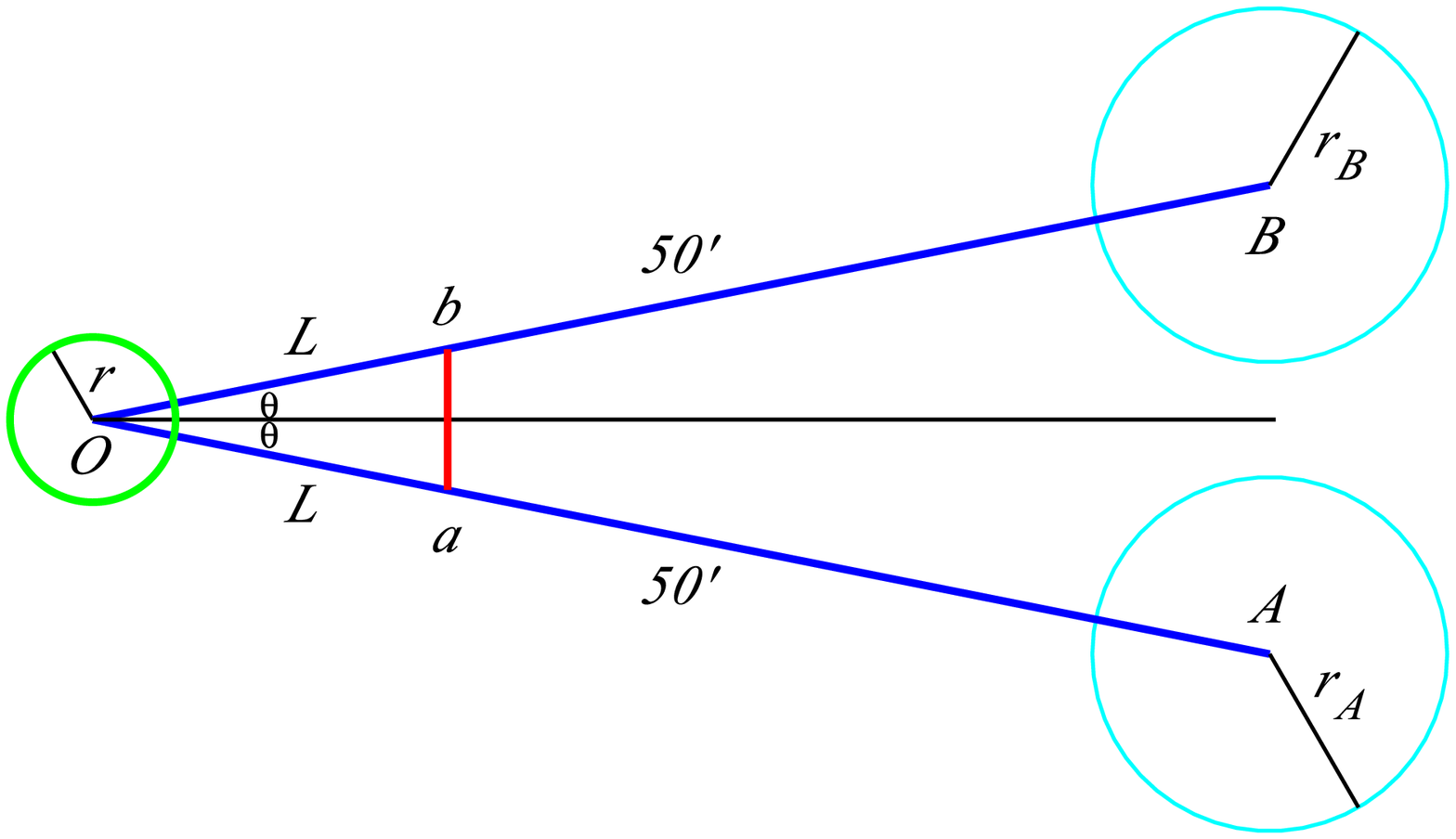}}
\caption{The digging site}
\end{figure}

Legrand states that points $a$ and $b$ are $2.5''=5/24'$ apart.  The
distances $aA$ and $bB$ are $50'$. 
Note that $\sin\theta=5/[48(r+L)]$.  The distance between the centres of
the holes is
\begin{equation}
AB=2(r+L+50)\sin\theta=\frac{5(r+L+50)}{24(r+L)}.\label{AB}
\end{equation}
The condition that the holes do not overlap is that $AB>r_A+r_B$.  
With \eqref{AB} this
reduces to 
\begin{equation}
r+L<\frac{250}{24(r_A+r_B)-5}.\label{nonovelap}
\end{equation}

We know from the story that $r_A=2'$.  The narrator says that
$r_B$ is slightly larger.  Even if we say $r_B=2'$,
for non-overlapping holes we must have $r+L<250/91\doteq 2'\,9''$.
This is clearly impossible.  The tree is 
\begin{quotation}
{\em
an enormously tall tulip tree, which stood, with some eight or ten oaks, upon the level, and far surpassed them all, and all other trees which I had then ever seen, in the beauty of its foliage and form, in the wide spread of its branches, and in the general majesty of its appearance.
}
\end{quotation}
And Poe states that Jupiter climbs $60$ or $70$ feet up to the branch with
the skull.  We should take $r$ to be the radius of the trunk at the height
of the skull branch.
For there to 
be a branch at that height substantial enough to hold a man's weight, 
$r$ must be at least $2$ inches and could easily be $6$ inches or more.
(A tulip tree on my campus was recently cut down.  The
stump measures $25''$ across.  It was not a particularly old tree.)  
But that means $L<2'\,7''$.  This puts the skull too close to
the tree trunk.  Jupiter crawls out to the end of the branch so this must
be more than $3$ feet.  Taking $r_B>2'$ only makes things worse.
Hence, although it's a great story, Poe neglected to do his trigonometry
homework and the two treasure
holes must have overlapped.  

From the wording in the story, it is not entirely clear whether the 
$50$ feet to the centre of the two holes is measured from the eye socket
drop points $a$ and $b$, as we have taken it, or from the tree trunk.
If from the trunk then the same assumptions on $r_A$, $r_B$ and $r$ as
above require $L$ be even smaller for the holes to not overlap.

\section{Background}
{\it The Gold-Bug} has been in the public domain for many years so there are
a great number of editions available, both print and electronic.
Excerpts here are taken from \cite{poe}.
For more on fiction that contains mathematics, see Alex Kasman's web page
\cite{kasman}.

In the story, Poe hints the treasure is Captain Kidd's.  William Kidd, born
circa 1654, was hanged for murder and piracy in England in 1701. He was 
certainly a privateer, which was legal in those times, and he held a
Letter of Marque signed by King William~III authorising him to attack 
French and certain other vessels.  It is uncertain
whether he was a pirate or not. After his death a legend grew that he
had left a treasure buried on the North American Atlantic coast.
For example, see \cite{harris}.

\end{document}